\documentclass{article}

\usepackage{amssymb,latexsym,amsmath}

\newtheorem{theorem}{Theorem}[section]
\newtheorem{lemma}[theorem]{Lemma}
\newtheorem{proposition}[theorem]{Proposition}

\title{The inverse problem for  linearly related orthogonal polynomials:  General case} 
\author{A. Pe\~{n}a $^{*}$, M. L. Rezola \thanks{Partially supported by MTM20015-65888-C4-4-P  (MINECO/FEDER) and DGA project E-64 (Spain).}\\ Departamento de Matem\'{a}ticas and IUMA \\ Universidad de Zaragoza (Spain)}

\begin{document}
\date{}

\maketitle

\bigskip

\begin{abstract}
We study the inverse problem in the theory of (standard) orthogonal polynomials involving two polynomials families $(P_n)_n$ and $(Q_n)_n$ which are connected by a linear algebraic structure such as 
$$P_n(x)+\sum_{i=1}^N r_{i,n}P_{n-i}(x)=Q_n(x)+\sum_{i=1}^M s_{i,n}Q_{n-i}(x)$$
for all $n=0,1, \dots$ where $N$ and $M$ are arbitrary nonnegative integer numbers.
\end{abstract}

\medskip


\medskip


\bigskip

\noindent {\it Address of the corresponding author}: \medskip
\\A. Pe\~{n}a\\Departamento de Matem\'{a}ticas\\Universidad de Zaragoza\\50009 Zaragoza (Spain)\\e.mail: anap@unizar.es.

\newpage

\section{Introduction} 

\setcounter{equation}{0}

The analysis of linear structure relations involving two monic orthogonal polynomial sequences (MOPS), $( P_n )_n$ and  $( Q_n )_n$,  such as 
\begin{equation}\label{relation N-M}
P_n(x)+\sum_{i=1}^N r_{i,n}P_{n-i}(x)=Q_n(x)+\sum_{i=1}^M s_{i,n}Q_{n-i}(x), \quad n \ge 0,
\end{equation}
where $N$ and $M$ are fixed nonnegative integer numbers, and $ (r_{i,n})_n$ and $(s_{i,n})_n$  are  sequences of complex numbers (and empty sum equals zero), has been a subject of research interest in the last decades.  In the literature, many works can be found where this type of relations is studied from different points of view. About the interest and importance of the study of these structure relations we refer to the introduction given in \cite{Petronilho} and \cite{APPR2013}, as well as the references therein.

In  many of these works the main problem stated and solved therein was the following inverse problem: assuming that $( P_n )_n$ is a MOPS and $( Q_n )_n$ only a simple set of polynomials ($Q_n$ is a polynomial of degree $n$),  verifying (\ref {relation N-M}), to find necessary and sufficient conditions so that $( Q_n )_n$ is also a MOPS and to obtain the relation between the corresponding regular linear functionals.  We want to notice that most of these papers deal with relations considering concrete values for $N$ and $M$ (see \cite{AMPR2003}, \cite{AMPR2010}, \cite{APRM2011}, \cite{APPR2013}, \cite{MP1995}, \cite{MV2014}). In this contribution we analyze the inverse problem for any values of  $N$ and  $M$. 

A classical tool for working with algebraic properties of orthogonal sequences of polynomials is the use of recurrence relations. Any MOPS $( P_n )_n$ is characterized by a three-term recurrence relation
\begin {equation*}
xP_n(x)=P_{n+1}(x)+\beta_nP_n(x)+\gamma_nP_{n-1}(x), \quad n=0,1,\dots
\end{equation*}
with initial conditions $P_{-1}(x)=0$ and $P_0(x)=1$, where $( \beta_n )_n$ and $( \gamma_n )_n$ are sequences of complex numbers such that $\gamma_n \not=0$ for all $n=1,2, \dots$.  This is known as Favard's theorem (see, e.g., \cite{Chihara}).

However this tool can be replaced by the use of dual basis which produces a natural way for studying the algebraic properties of sequences of orthogonal polynomials. Any simple set of polynomials $( P_n )_n$  has a dual basis $({\bf a}_n)_n$,  that is 
$$\langle {\bf a}_n, P_j \rangle:=\delta_{n,j},\quad  n,j=0,1, \dots$$
being $\delta_{n,j}$ the usual Kronecker symbol. Moreover, if $( P_n )_n$ is a MOPS with respect to the linear functional ${\bf u}$, then the associated dual basis is 
$${\bf a}_n=\frac{P_n}{\langle {\bf u}, P_n^2 \rangle}\,{\bf u}, \quad n=0,1, \dots$$
see \cite{Maroni1991}.  A more detailed description can be seen in Section 2 of \cite{Petronilho}.

This is the main tool used by Petronilho in  \cite{Petronilho} to solve part of the inverse problem for general relations as (\ref {relation N-M}) assuming the orthogonality of both sequences $( P_n )_n$ and $( Q_n )_n$ with respect to ${\bf u}$ and ${\bf v}$, as well as some additional assumptions which we will call initial conditions.  More precisely, he obtains that the rational transformation between the linear functionals ${\bf u}$ and ${\bf v}$ is $${\Phi}_M {\bf u}={\Psi}_N{ \bf v}$$ where ${\Phi}_M$ and ${\Psi}_N$ are polynomials of degrees $M$ and $N$, respectively. In Section 2, we make an exhaustive study about these initial conditions.  We obtain, in Theorem \ref{caracterizacion condiciones Petro}, that these initial conditions characterize the existence of such rational transformation and the non existence of another relationship with degrees less than or equal to $M$ and $N$, respectively. Moreover these initial conditions allow us to prove important facts like:

 (a) All the coefficients $ r_{N,n}$ and $ s_{M ,n}$ in (\ref {relation N-M}) are not zero for $n \ge N+M$.
 
 (b)  There exist  constant sequences whose values are the coefficients of the polynomials ${\Phi}_M$ and ${\Psi}_N$.  It is important to note that the existence of such constant sequences was already obtained in \cite{MP1995} for $N=1, M=0$, in \cite{AMPR2003} for $N=1,M=1$, in \cite{APRM2011} for $N=2,M=0$ and in \cite{APPR2013} for $N=2,M=1$, in terms of the recurrence coefficients.  Here, for arbitrary values of $N$ and $M$, we obtain this property in a compact form by using linear functionals and determinants. 
 
 \medskip

On the other hand, in Section 3 we get necessary and sufficient conditions in order to  the sequence $( Q_n )_n$ defined recursively by  (\ref {relation N-M}) becomes also a MOPS. The main advance in this line is to introduce some auxiliary polynomials $R_n$, which are precisely the linear combinations of the polynomials $P_n$ that appear in the relation (\ref {relation N-M}). These polynomials do not necessarily have to be orthogonal, but they are interesting in two senses. First, they allow to simplify the computations in the problem of characterizing the orthogonality of the sequence $( Q_n )_n$. And second and more important, is that the conditions which characterize the orthogonality of $( Q_n )_n$  correspond to most of the conditions ($n \ge N+M+1$) that characterize the orthogonality of two simpler problems. More precisely, assuming that $( P_n )_n$ is a MOPS, to characterize the orthogonality of $( R_n )_n$ and then, assuming that $( R_n )_n$ is a MOPS,  to characterize the orthogonality of $( Q_n )_n$. In some way the problem  $N-M$ can be divided in two simpler problems $N-0$ and $0-M$ but always keeping in mind that not all the conditions of regularity (orthogonality) appear, the first ones are different. This is because, as we have already mentioned, these auxiliary polynomials $( R_n )_n$ do not have to be orthogonal.

\section{Relation between the regular functionals and consequences} 

\setcounter{equation}{0}

 Let   $( P_n )_n$ and  $( Q_n )_n$ be two sequences of monic polynomials  orthogonal with respect to the regular functionals $\bf u$ and $\bf v$, respectively, normalized by $\langle {\bf u}, 1 \rangle=1= \langle {\bf v}, 1 \rangle$.  Suppose that these families of polynomials are  linearly related by (\ref{relation N-M}), that is
$$P_n(x)+\sum_{i=1}^N r_{i,n}P_{n-i}(x)=Q_n(x)+\sum_{i=1}^M s_{i,n}Q_{n-i}(x), \quad n \ge 0.$$
Consider the following auxiliary polynomials, namely $R_n$, 
$$R_n(x):=P_n(x)+\sum_{i=1}^N r_{i,n}P_{n-i}(x)=Q_n(x)+\sum_{i=1}^M s_{i,n}Q_{n-i}(x). $$
Denote by $({\bf c }_n)_n$, $({\bf a}_n)_n$ and $({\bf b}_n)_n$ the dual basis corresponding to  $( R_n )_n$, $( P_n )_n$ and  $( Q_n )_n$, respectively. Expanding $({\bf a}_n)_{n=0}^{M-1}$ and $({\bf b}_n)_{n=0}^{N-1}$ in terms of $({\bf c}_n)_n$, we can write
\begin{equation} \label{def A}
{\bf A} ({\bf c }_0, \dots, {\bf c }_{M+N-1})^T=({\bf a }_0, \dots, {\bf a }_{M-1}, {\bf b }_0, \dots, {\bf b }_{N-1})^T
\end{equation}
where ${\bf A}$ is a $(M+N) \times (M+N)$ matrix whose elements are the coefficients that appear in the relation (\ref{relation N-M}). Its explicit expression can be seen in \cite[Theorem 1.1]{Petronilho}. There, it was proved that the initial conditions 
\begin{equation}\label{initial conditions}
\text{det} \,{\bf A} \not=0, \quad r_{N,N+M}\not=0, \quad \rm{and} \quad s_{M,N+M}\not=0 
\end{equation}
yield the following  relation between the linear functionals  $\bf u$ and $\bf v$ 
\begin{equation}\label{relacion funcionales}
{\Phi}_M {\bf u}={\Psi}_N{ \bf v}
\end{equation}
where ${\Phi}_M$ and ${\Psi}_N$ are polynomials of (exact) degrees $M$ and $N$, respectively. 

\medskip

In the sequel we will use the following notations:
$$\overline{P}_n(x)=\frac{P_n(x)}{\langle {\bf u}, {P_n}^2 \rangle}, \quad \text  {and} \quad
\overline{Q}_n(x)=\frac{Q_n(x)}{\langle {\bf v}, {Q_n}^2 \rangle}.$$

For $n\ge N-1$,  we introduce the $N \times N$  matrices $\bf B_n$ and ${\bf B}_{n}^{ i},\,i=0, \dots, N-1,$ where

\begin{equation*}\label{matrix Bn}
{\bf B}_n:= \begin{pmatrix}
 \langle \overline{Q}_{N-1}{\bf v}, P_{n}\rangle & \cdots & \langle \overline{Q}_0{\bf v}, P_{n}\rangle \\
 \langle \overline{Q}_{N-1}{\bf v}, P_{n-1}\rangle  &\cdots & \langle \overline{Q}_0{\bf v}, P_{n-1}\rangle \\
\vdots &  \ddots & \vdots \\
\langle \overline{Q}_{N-1}{\bf v}, P_{n-(N-1)}\rangle &\cdots & \langle \overline{Q}_0{\bf v}, P_{n-(N-1)}\rangle \\
\end{pmatrix}_{N \times N} 
 \end{equation*}
and ${\bf B}_{n}^{ i}:= $ the matrix  obtained replacing  in  ${\bf B}_{n}$  the $ i$-th column by  the vector $ (\langle \overline{Q}_{N}{\bf v}, P_{n}\rangle, \dots,  \langle \overline{Q}_{N}{\bf v}, P_{n-(N-1)}\rangle)^T$.

In a similar way for  $n\ge M-1$,  we introduce the $M \times M$  matrices $ \widetilde {\bf B}_n$ and 
${\widetilde{\bf B}}_n^{i}, \, i=0, \dots, M-1,$ where 

\begin{equation*}\label{matrix tildeBn}
 \widetilde {\bf B}_n:=\begin{pmatrix}
 \langle \overline{P}_{M-1}{\bf u}, Q_{n}\rangle & \cdots & \langle \overline{P}_0{\bf u}, Q_{n}\rangle \\
 \langle \overline{P}_{M-1}{\bf u}, Q_{n-1}\rangle  &\cdots & \langle \overline{P}_0{\bf u}, Q_{n-1}\rangle \\
\vdots &  \ddots & \vdots \\
\langle \overline{P}_{M-1}{\bf u}, Q_{n-(M-1)}\rangle &\cdots & \langle \overline{P}_0{\bf u}, Q_{n-(M-1)}\rangle \\
\end{pmatrix}_{M \times M} 
 \end{equation*}
and ${\widetilde{\bf B}}_n^{i}:=$  the matrix  obtained replacing  in  $\widetilde{\bf B}_n$ the $i$-th column by the vector 
 $ \left( \langle \overline{P}_{M}{\bf u}, Q_{n}\rangle, \dots,  \langle \overline{P}_{M}{\bf u}, Q_{n-(M-1)}\rangle \right)^T.$ 
 
 \medskip

Next, we show a property about  the determinants of  the matrices $\bf B_n$ and $ \widetilde {\bf B}_n$ which plays an important role in this work.

  \begin{lemma}\label{relation entre Bn} Let   $( P_n )_n$ and  $( Q_n )_n$ be two sequences of monic polynomials linearly related by  (\ref {relation N-M}). 
    \item{(a)} If the polynomials $Q_n$ are orthogonal with respect to the functional {\bf v},  then
  $$\text{det} \,\, {\bf B}_{n}=(-1)^N\, r_{N,n}\,\text{det}\,\,{\bf B}_{n-1}, \quad {\text for} \quad  n \ge M+N.$$
 
  \item{(b)} If the polynomials $P_n$ are orthogonal with respect to the functional {\bf u},  then
   $$\text{det} \,\,{\widetilde{\bf B}}_{n}=(-1)^M\, s_{M,n}\,\text{det}\,\,{\widetilde{\bf B}}_{n-1}, \ \quad {\text for} \quad  n \ge M+N.$$
\end{lemma}
  \textbf{Proof.} (a) First we get that $\text{det} \,{\bf B}_{N+M}=(-1)^N\, r_{N,N+M}\,\text{det}\,{\bf B}_{N+M-1}$. Indeed, taking into account  that     
  $$P_{N+M}(x)=Q_{N+M}(x)+\sum_{i=1}^M s_{i,N+M}Q_{N+M-i}(x)-\sum_{i=1}^N r_{i,N+M}P_{N+M-i}(x)$$  by (\ref {relation N-M}) for $n=N+M$, 
using the orthogonality of the polynomials $Q_n$ with respect to  {\bf v} and  developing the determinant of ${\bf B}_{N+M}$ by the first row, it can be derived that
   $$ \text{det} \,{\bf B}_{N+M}= - \,(-1)^{N-1}\, r_{N,N+M} \, \text{det} \,{\bf B}_{N+M-1}.$$
   To conclude the proof of $(a)$  it suffices to observe that the same argument works for any fixed $n >N+M$. 
   
   (b) This property can be similarly derived with the appropriate changes.       \hfill$\Box$
  
   \bigskip
   
  Now, we will see that the polynomials $\Psi_N$  and $\Phi_M$ which satisfy the relation  (\ref {relacion funcionales}) can be written as:
 
\begin{equation}\label{def Psi}
\Psi_N(x)=r_{N,M+N}\,\overline{Q}_N(x)+ \sum_{i=0}^{N-1}  \lambda_i\, \overline{Q}_i(x),
\end{equation}
and 
\begin{equation}\label{def Phi}
\Phi_M(x)=s_{M,M+N} \,\overline{P}_M(x)+ \sum_{i=0}^{M-1} \mu_i\, \overline{P}_i(x).
\end{equation}
Indeed,  if $\Psi_N$ is a polynomial of degree $N$ it can be written as 
$$\Psi_N(x)=\lambda_N\,\overline{Q}_N(x)+ \lambda_{N-1}\, \overline{Q}_{N-1}(x)+\dots+\lambda_1\, \overline{Q}_1+\lambda_0$$ with $\lambda_N\not=0$. Then, using the dual basis of  $( P_n )_n$  we have 
\begin{equation}\label{base dual (P_n)}
\Psi_N {\bf v}=\sum_ {j=0}^M \langle {\Psi}_N {\bf v}, P_j \rangle\, {\overline  P_j} {\bf u}
\end{equation} 
because  $\langle {\Psi}_N {\bf v}, P_j \rangle=\langle {\Phi}_M {\bf u}, P_j \rangle=0 \quad \text {for} \quad  j\ge M+1$.  Moreover, by (\ref{relation N-M}) $r_{N,N+M}\, \langle {\Psi}_N {\bf v}, P_M \rangle=s_{M,N+M}\, \langle {\Psi}_N {\bf v}, Q_N \rangle=\lambda_N \,s_{M,N+M}$ and from  (\ref {base dual (P_n)})  we get the expressions (\ref{def Psi}) and (\ref{def Phi}) for the polynomials ${\Psi}_N$ and ${\Phi}_M$, respectively.    

\medskip

Next, we obtain a characterization of the initial conditions (\ref{initial conditions}) in terms of the relations between the linear functionals ${\bf u}$ and ${\bf v}$.

 \begin{theorem}\label{caracterizacion condiciones Petro}
 Let   $( P_n )_n$ and  $( Q_n )_n$ be two MOPSs  with respect to the regular functionals $\bf u$ and $\bf v$, respectively, normalized by $\langle {\bf u}, 1 \rangle=1= \langle {\bf v}, 1 \rangle$.  Assume that they are  related by (\ref{relation N-M}) where $N,M \ge 1$.
 
 Then, the following statements  are equivalent:
 
\noindent{(a)} $$ \text{det} \,{\bf A} \not=0, \quad r_{N,N+M}\not=0, \quad \rm{and} \quad s_{M,N+M}\not=0$$
(b) There exist polynomials  $\Phi_M, \Psi_N$ of degrees $M, N$, respectively, such that  ${\Phi}_M\,{\bf u}= \Psi_N\,{\bf v}$, and there is no other relationship $\widetilde{\Phi}_M\,{\bf u}= \widetilde{\Psi}_N\,{\bf v}$ with degrees of polynomials $\widetilde{\Phi}_M$ and $\widetilde{\Psi}_N$ less than or equal to 
$M$ and $N$ respectively.
 \end{theorem}
 
  \textbf{Proof.} $(a) \Longrightarrow (b)$
  
  First we prove that the polynomials $\Psi_N$  and $\Phi_M$  defined by  (\ref{def Psi}) and  (\ref{def Phi}) such that ${\Phi}_M\,{\bf u}= \Psi_N\,{\bf v}$ are unique. Indeed, using the dual basis $({\bf a}_n)_n$ and $({\bf b}_n)_n$ corresponding to $(P_n)_n$ and  $( Q_n )_n$, respectively, we have 
  
\begin{equation*}\label{def Psi-v}
\Psi_N\,{\bf v}=r_{N,M+N}\,{\bf b}_N+ \lambda_{N-1}\, {\bf b}_{N-1}+\dots+\lambda_1\,{\bf b}_1+\lambda_0\,{\bf b}_0,
\end{equation*}
and
 \begin{equation*}\label{def Phi-u}
{\Phi}_M\,{\bf u}=s_{M,M+N} \,{\bf a}_M+ \mu_{M-1}\, {\bf a}_{M-1}+\dots+\mu_1\, {\bf a}_1+\mu_0\,{\bf a}_0,
\end{equation*}
since ${\bf a}_n=\overline{P}_n\,{\bf u}$ and  ${\bf b}_n=\overline{Q}_n\,{\bf v}$ for $n \ge 0.$ 

 So, by  the relation between the functionals (\ref{relacion funcionales})  and the definition (\ref{def A}) of the matrix ${\bf A}$, we get
   \begin{align*}
   &r_{N,M+N}\,{\bf b}_N-s_{M,M+N} \,{\bf a}_M=\\
   &(\mu_0, \dots,  \mu_{M-1}, -\lambda_0, \dots, -\lambda_{N-1})({\bf a}_0, \dots, {\bf a}_{M-1},{\bf b}_0, \dots, {\bf b}_{N-1})^T= \\
   & (\mu_0, \dots,  \mu_{M-1}, -\lambda_0, \dots, -\lambda_{N-1})\,{\bf A}\,({\bf c }_0, \dots, {\bf c }_{M+N-1})^T. \\
   \end{align*} 
 Thus, since  $({\bf c}_n)_n$ is a basis and  there exists ${\bf A}^{-1}$ the inverse matrix of ${\bf A}$, the vector $(\mu_0, \dots,  \mu_{M-1}, -\lambda_0, \dots, -\lambda_{N-1})$ is unique. Therefore there exist unique polynomials $\Psi_N$  and $\Phi_M$  of degrees respectively $N,M$ such that ${\Phi}_M\,{\bf u}= \Psi_N\,{\bf v}$. 
 
 Now,  we suppose that there exist another polynomials $\widetilde{\Psi}_N$  and $\widetilde{\Phi}_M$  satisfying  $\widetilde{\Phi}_M\,{\bf u}= \widetilde{\Psi}_N\,{\bf v}$ with $\text {deg}\,\,\widetilde{\Psi}_N \le N$ and   $\text {deg}\,\,\widetilde{\Phi}_M \le M.$ Thus, we have
 $$({\Phi}_M - \widetilde{\Phi}_M)\,{\bf u}=({\Psi}_N - \widetilde{\Psi}_N)\,{\bf v}$$
 so it is not possible  that  $\text {deg}\,\,\widetilde{\Psi}_N<N$  and  $\text{deg}\,\,\widetilde{\Phi}_M<M$ hold simultaneously. 

To conclude,  it suffices to observe that from the two relations between the functionals ${\bf u}$ and ${\bf v}$, we also obtain the  following relation  
 $$\widetilde{\Psi}_N\,{\Phi}_M\,{\bf v}=\widetilde{\Phi}_M\,\Psi_N\,{\bf v}$$ 
 which yields a contradiction if  either $\text{deg}\,\,\widetilde{\Phi}_M=M$ and $\text{deg}\,\,\widetilde{\Psi}_N<N$ or $\text{deg}\,\,\widetilde{\Phi}_M<M$ and $ \text{deg}\,\,\widetilde{\Psi}_N=N.$
 
  $(b) \Longrightarrow (a)$
  
  Consider the systems 
  \begin{equation}\label {sistema Psi}
 \langle {\Psi}_N {\bf v}, P_n \rangle=0, \quad n= M+1, \dots, N+M 
 \end{equation}
  and 
 \begin{equation*}\label {sistema Phi}
 \langle {\Phi}_M {\bf u}, Q_n \rangle=0, \quad n= N+1, \dots, N+M 
\end{equation*}
 where the unknowns are, respectively, $(\lambda_i)_{i=0}^{N-1}$ and $(\mu_i)_{i=0}^{M-1}.$
 Notice that, by hypothesis,  the solutions of these systems are unique. Thus, the respective  coefficient matrices  ${\bf B}_{N+M}$  and  $\widetilde{\bf B}_{N+M}$ have maximum rank that is
 $${\text det} \,\,{\bf B}_{N+M} \not=0, \quad \text {and} \quad {\text det} \,\,\widetilde{\bf B}_{N+M} \not=0.$$ So,  by Lemma  \ref{relation entre Bn} we obtain $r_{N,M+N}\not=0$ and $s_{M,M+N}\not=0$.
 
 Now, it remains only to prove that $\text{det} \,{\bf A} \not=0.$  To do this, we consider the system with two equations, one of which is the expansion of ${\bf b}_N$ as a linear combination of $({\bf c}_i)_{i=0}^{M+N}$ and the other one is the corresponding expansion of ${\bf a}_M$. Then, multiplying the first of these equations by $r_{N,M+N}$ and the second one by $s_{M,M+N}$, and subtracting the resulting equations (this will eliminate ${\bf c }_{M+N}$) we get
 $$r_{N,M+N}\,{\bf b}_N-s_{M,M+N} \,{\bf a}_M=\sum_{i=0}^{M+N-1} X_i \,{\bf c}_i$$
 (see the proof of Theorem 1.1 in \cite{Petronilho}, for instance). On the other hand, we know 
$$r_{N,M+N}\,{\bf b}_N-s_{M,M+N} \,{\bf a}_M=(\mu_0, \dots,  \mu_{M-1}, -\lambda_0, \dots, -\lambda_{N-1})\,{\bf A}\,({\bf c }_0, \dots, {\bf c }_{M+N-1})^T.$$ Thus, 
 $${\bf A}^T \,(\mu_0, \dots,  \mu_{M-1}, -\lambda_0, \dots, -\lambda_{N-1})^T=( X_0, \dots X_{M+N-1})^T$$ and  then $\text{det} \,{\bf A}=\text{det} \,{\bf A}^T \not=0$
because the system has a unique solution.   \hfill$\Box$

  \medskip
 
In the following theorem we prove that the initial conditions (\ref{initial conditions}) allow us to assure that the lengths of the linear combinations of $( P_n )_n$ and  $( Q_n )_n$ in (\ref{relation N-M})   are exactly $N+1$ and $M+1$ respectively.  Moreover, we can find constant sequences whose values are precisely the coefficients of the polynomials  $\Psi_N$  and $\Phi_M$.

 \begin{theorem}\label{unicidad Psi-Phi}
 Let   $( P_n )_n$ and  $( Q_n )_n$ be two MOPSs  with respect to the regular functionals $\bf u$ and $\bf v$, respectively, normalized by $\langle {\bf u}, 1 \rangle=1= \langle {\bf v}, 1 \rangle$. Assume that they are  related by (\ref{relation N-M})  where $N,M \ge 1$ and the coefficients satisfy the initial conditions (\ref{initial conditions}). 
  Then, the following properties hold:
\item{(a)} All the coefficients  $r_{N,n}$ and $s_{M,n}$ in (\ref{relation N-M}) are not zero, for every  $n \ge N+M$. 
\item{(b)} There exist constant sequences $(\lambda_{i,n})_{n\ge N+M}$ and $(\mu_{i,n})_{n\ge N+M}$ such that 
$\lambda_{i,n}=\lambda_i,\, i=0, \dots, N-1$ and $\mu_{i,n}=\mu_i,\,  i=0, \dots, M-1$ where $(\lambda_i)_{i=0}^{N-1}$ and $(\mu_i)_{i=0}^{M-1}$ are the coefficients which appear in the expressions (\ref{def Psi}) and (\ref{def Phi}) of the polynomials  $\Psi_N$  and $\Phi_M$.
\end{theorem}

 \textbf{Proof.} (a) First we observe that for $n \ge N+M$ we have
 \begin{align*}
 &\frac{s_{M,N+M}}   {\langle {\bf u}, P_M^2 \rangle}\,r_{N,n}
 \, \langle {\bf u}, P_{n-N}^2 \rangle =
 \langle {\Phi}_M {\bf u},  (P_n+\sum_{i=1}^N r_{i,n}P_{n-i})\,Q_{n-(N+M)} \rangle \\
 &=\langle  {\Psi}_N{ \bf v}, (Q_n+\sum_{i=1}^M s_{i,n}Q_{n-i})\,Q_{n-(N+M)} \rangle
 =\frac{r_{N,N+M}} {\langle {\bf v}, Q_N^2 \rangle}\,s_{M,n} \, \langle {\bf v}, Q_{n-M}^2 \rangle
 \end{align*}
 taking into account the hypothesis  (\ref{relation N-M}), (\ref{relacion funcionales}), the expressions (\ref{def Psi})-(\ref {def Phi}) and the orthogonality of the polynomials $P_n$ and $Q_n$ with respect to the functionals ${\bf u}$ and ${\bf v}$, respectively. Thus, since $r_{N,N+M}\,s_{M,N+M} \not=0$,   it is enough to prove that either $r_{N,n}\not=0$ or $s_{M,n}\not=0$ for  $n \ge N+M+1$. 
 
 Here, we will prove that $r_{N,n}\not=0$ for $n \ge N+M+1$. 

 Assume that there exists $n \ge N+M+1$ such that $r_{N,n}=0$ and let $n_0:=\min \{n; n\ge N+M+1, r_{N,n}=0\}$. Then Lemma \ref {relation entre Bn} gives $\text {det} \, {\bf B}_{n_0}=0$, so the coefficient matrix associated with the system 
 $$\langle {\Psi}_N {\bf v}, P_n \rangle=0, \quad  n=n_0-(N-1), \dots, n_0$$ has not maximum rank and therefore there exists another solution, namely, $(\tilde{\lambda}_i)_{i=0}^{N-1}$  such that 
 $$\langle  \widetilde{\Psi}_N {\bf v}, P_n \rangle=0, \quad  n= n_0-(N-1), \dots, n_0$$
 where 
 $$ \widetilde{\Psi}_N (x)=r_{N,M+N}\,\overline{Q}_N(x)+\tilde{\lambda}_{N-1}\, \overline{Q}_{N-1}(x)+\dots+\tilde{\lambda}_1\, \overline{Q}_1+\tilde{\lambda}_0.$$ 
 Moreover  $\langle  \widetilde{\Psi}_N {\bf v}, P_{n_0-N} \rangle\not=0$. Indeed, if $\langle  \widetilde{\Psi}_N {\bf v}, P_{n_0-N} \rangle=0$ then the system 
 $$\langle  \widetilde{\Psi}_N {\bf v}, P_n \rangle=0, \quad  n= n_0-N, \dots, n_0-1$$ has two solutions which yields a contradiction since the coefficient matrix of this system, ${\bf B}_{n_0-1}$, has  maximum rank by definition of $n_0$.
 
 Hence, using the relation (\ref{relation N-M}), we have
 $$\langle  (\Psi_N-\widetilde{\Psi}_N) {\bf v}, P_{n} \rangle= -\langle  \widetilde{\Psi}_N {\bf v}, P_{n} \rangle=0, \quad n \ge n_0-(N-1),$$ 
 and   $$\langle  (\Psi_N-\widetilde{\Psi}_N) {\bf v}, P_{n_0-N} \rangle= -\langle  \widetilde{\Psi}_N {\bf v}, P_{n_0-N} \rangle\not=0.$$    
 Denote by   $h_{N-1}$  the polynomial $\Psi_N-\widetilde{\Psi}_N$ of degree less than or equal to $N-1$. Then, writting   $h_{N-1}{\bf v}$ in the dual basis of  $( P_n )_n$,  we have 
  $$h_{N-1}{\bf v}=\sum_{j=0}^{n_0-N} \frac{\langle h_{N-1} {\bf v}, P_j \rangle}{\langle {\bf u}, P_j^2 \rangle}P_j {\bf u},$$ so there exists a polynomial $\varphi_{n_0-N}$ of degree $n_0-N$ such that  $h_{N-1} {\bf v}=\varphi_{n_0-N} {\bf u}$.  
  
 Finally observe that  since the functionals ${\bf u}$ and ${\bf v}$ satisfy the  two relations
 $$h_{N-1} {\bf v}=\varphi_{n_0-N} \,{\bf u} \quad \text {and} \quad {\Psi}_N{ \bf v}={\Phi}_M {\bf u} ,$$ it  can be obtained
 $$h_{N-1} {\Phi}_M {\bf v}=\varphi_{n_0-N} {\Psi}_N{ \bf v}$$ which  leads to a contradiction,  taking into account the degrees of the polynomials $h_{N-1} {\Phi}_M$ and $ \varphi_{n_0-N} {\Psi}_N$.

 (b)  From the previous theorem, we already know that the system  (\ref{sistema Psi})
 $$\langle {\Psi}_N {\bf v},  P_n\rangle=0, \, n=M+1, \dots, N+M$$ 
 has a unique solution, namely $(\lambda_i)_{i=0}^{N-1}$, which are the coefficients of the polynomial $\Psi_N$.
 So, since  $ \text{det} \,{\bf B}_{N+M}\not=0$ by Cramer's rule 
 $$\lambda_i=-r_{N,N+M}\, \frac{\text{det} \,\,{\bf B}_{N+M}^{ i}}{\text{det} \,\, {\bf B}_{N+M}}, \quad \text {for} \quad  i=0,\dots,N-1.$$ 
 On the other hand, since $ {\Phi}_M {\bf u}= {\Psi}_N {\bf v}$, it is obvious that 
  $$\langle {\Psi}_N {\bf v},P_n\rangle=0, \, n \ge M+1.$$
 Now, for each $n$ fixed,  $n \ge N+M$, we can consider the system
\begin{equation*}\label{sistema Psi n grande}
\langle {\Psi}_N {\bf v}, P_i\rangle=0, \, i=n-(N-1), \dots, n
\end{equation*} 
whose matrix of coefficients  ${\bf B}_n$ has  maximum rank, because we know   
$$\text{det} \,\,{\bf B}_n=(-1)^N\,r_{N,n}\dots \,
(-1)^N\,r_{N,N+M+1}\text{det} \,\,{\bf B}_{N+M}\not=0,$$ by (a) and  Lemma  \ref{relation entre Bn}. Then, for every $n \ge N+M$, this system has a unique solution for $(\lambda_i)_{i=0}^{N-1}$, namely  $(\lambda_{i,n})_{i=0}^{N-1}$, and again by Cramer's rule we have$$\lambda_{i,n}=-r_{N,N+M}\, \frac{\text{det} \,\,{\bf B}_{n}^{ i}}{\text{det} \,\, {\bf B}_{n}}, \quad \text {for} \quad  i=0,\dots,N-1,$$
for all $n\ge N+M.$

Thus,  we can assure there exist  $N$ constant sequences, namely $(\lambda_{i,n})_{i=0}^{N-1}$,  such that  for every $n \ge N+M$ we have $\lambda_{i,n}=\lambda_i$  that is the value of these constants coincide with the coefficients of the polynomial $\Psi_N$.

To conclude the proof of $(b)$ , we  work in the same way  with the polynomial $\Phi_M$.  So, we have 
$$\langle \Phi_M {\bf u}, Q_n \rangle=0, \, n\ge N+1.$$ Thus, for each n fixed,  $n \ge N+M$ we can consider the system
\begin{equation*}\label{sistema Phi n grande}
\langle {\Phi}_M {\bf u},  Q_i\rangle=0, \, i=n-(M-1), \dots, n
\end{equation*}
where the unknowns are $(\mu_i)_{i=0}^{M-1}$,  that is the coefficients of the polynomial ${\Phi}_M$, see (\ref {def Phi}).  The uniqueness of the polynomial ${\Phi}_M$ obtained in the previous theorem leads us to state that $\text{det} \,\, \widetilde{\bf B}_n\not=0,$ for all  $n \ge N+M$. Then, in the same way as before, we obtain  that there exist  $M$ constant sequences, namely $(\mu_{i,n})_{i=0}^{M-1}$,  such that  for every $n \ge N+M$ we have $\mu_{i,n}=\mu_i$  where 
$$\mu_{i,n}=-s_{M,N+M}\, \frac{\text{det} \,\,\widetilde{\bf B}_{n}^{ i}}{\text{det} \,\, \widetilde{\bf B}_{n}}.$$ 
Note that the values of these constants coincide with the coefficients of the polynomial $\Phi_M$.
  \hfill$\Box$

\medskip  

Observe that the notation used before Lemma  2.1 does not work for either  $N=0$ or $M=0$. So we conclude this section, for the sake of completeness, showing the analogous results  for these particular situations. We will show only the case $M = 0$ and $N \ge 1$, that is

\begin{equation}\label{relation N-0}
P_n(x)+\sum_{i=1}^N r_{i,n}P_{n-i}(x)=Q_n(x) \quad n \ge 0,
\end{equation}
where $N \ge 1$ and $ r_{i,n}$ ($i=1, \dots, N)$ are  complex numbers. The other case is totally  analogous with the appropriate changes.

\begin{theorem} \label{N-0}
 Let   $( P_n )_n$ and  $( Q_n )_n$ be two sequences of monic polynomials  with respect to the regular functionals $\bf u$ and $\bf v$, respectively, normalized by $\langle {\bf u}, 1 \rangle=1= \langle {\bf v}, 1 \rangle$.  Suppose that these families of polynomials are  linearly related by (\ref{relation N-0})  with the initial condition $r_{N,N}\not=0$. Then, the following properties hold:
 
 (a) The exists a unique polynomial $\Psi_N$ of degree $N$ such that ${\bf u}={\Psi}_N\,{\bf v}$ where  ${\Psi}_N$ is defined by (\ref{def Psi}) with $M=0$. 
 
 (b) All the coefficients $r_{N,n}$  in (\ref{relation N-0}) are not zero, for every $n \ge N.$
 
 (c) There exist constant sequences $(\lambda_{i,n})_{n\ge N}$  such that 
$\lambda_{i,n}=\lambda_i,\, i=0, \dots, N-1$  where $(\lambda_i)_{i=0}^{N-1}$ are the coefficients which appear in the expression (\ref{def Psi})  of the polynomial  $\Psi_N$.  
 \end{theorem}
 
  \textbf{Proof.}  The proof follows the same ideas of the previous theorems. In any case, for the sake of completeness, we briefly expose the arguments used. 
  
  Observe that
$\langle {\bf u}, Q_n \rangle=\langle {\bf u}, P_n+\sum_{i=1}^N r_{i,n}P_{n-i} \rangle=0$ for $n \ge N+1$. Then, using the dual basis of $( Q_n )_n$ we get
 ${\bf u}=\sum _{n=0}^N \, \langle {\bf u}, Q_n \rangle\, \overline{Q}_n\,{\bf v}$ 
 where $\langle {\bf u}, Q_N \rangle=r_{N,N}\not=0$ and therefore there exists a unique polynomial of degree $N$ such that $${\bf u}={\Psi}_N\,{\bf v}$$ where  ${\Psi}_N$ is defined by (\ref{def Psi}).  
 
On the other hand, we can check that ${\text det} \,\,{\bf B}_{N-1}=1$ and that $(a)$ of the Lemma \ref {relation entre Bn} holds with $M=0$. Thus, the coefficient matrix  ${\bf B}_N$ of the system $$\langle {\Psi}_N\,{\bf v}, P_n \rangle=0, \, n=1, \dots, N$$
has maximum rank  and  the coefficients $(\lambda_i)_{i=0}^{N-1}$ of the polynomial  $\Psi_N$ are  
$$\lambda_i=-r_{N,N}\, \frac{\text{det} \,\,{\bf B}_{N}^{ i}}{\text{det} \,\, {\bf B}_{N}}, \quad \text {for} \quad  i=0,\dots,N-1.$$ 
Moreover,  in the same way as in the previous theorem,
it can be proved that  $$r_{N,N}\not=0 \Longrightarrow  r_{N,n}\not=0, \quad  n\ge N.$$ 
Thus,  $\text{det} \,\, {\bf B}_{n} \not=0$ for $n\ge N$ and $$\lambda_i=\lambda_{i,n}=-r_{N,N}\, \frac{\text{det} \,\,{\bf B}_{n}^{ i}}{\text{det} \,\, {\bf B}_{n}}, \quad \text {for} \quad  i=0,\dots,N-1, \quad n\ge N.$$  \hfill$\Box$

\section{Orthogonality  characterizations}
\setcounter{equation}{0}

  Let   $( P_n )_n$ and  $( Q_n )_n$ be two sequences of monic polynomials linked by a structure relation  as (\ref{relation N-M}),  with the conventions $r_{N,n}s_{M,n}\not=0$ for all $n \ge N+M$.
  
  In this section, we want to find necessary and sufficient conditions in order to  $( Q_n )_n$ be  a MOPS if $( P_n )_n$ is a MOPS.  As we have mentioned in the introduction, to get this general case, we introduce some auxiliary polynomials $R_n$ and so, in some way, the problem $N-M$ can be divided in two simpler problems $N-0$ and $0-M$. So, we previously study these particular situations that correspond to consider in the relation (1.1) either $M=0$ or $N=0$.

From now on, $( P_n )_n$ denotes a MOPS with respect to a regular functional $\bf u$ and  $(\beta_n)_n$ and $(\gamma_n)_n$  the corresponding sequences of recurrence coefficients,  that is

\begin{align} \label{recurrenciaPn}
P_{n+1}(x)&=(x-\beta_n)P_n(x)-\gamma_nP_{n-1}(x), \quad n\ge0,
\\P_0(x)&=1, \quad P_{-1}(x)=0,\nonumber
\end{align}
\noindent with  $\gamma_n \not=0$ for all $n\ge1$.

\begin{proposition}\label{caracterizacion ortogonalidad Rn} Let $( P_n )_n$ be a MOPS with recurrence
coefficients $(\beta_n)_n$ and $(\gamma_n)_n$.  Fixed $N \ge 1$, we define a sequence
$(R_n)_n$ of monic polynomials by 
\begin{equation} \label{relation Rn-Pn}
R_n(x):=P_n(x)+\sum_{i=1}^N r_{i,n}P_{n-i}(x)
\quad n\ge0,
\end{equation}
\noindent where  $ r_{i,n}$ are complex numbers  such that $r_{i,n}=0$ if $i >n$  and $r_{N,n} \not=0$ for all $n \ge N$. 

Then $( R_n )_n$  is a MOPS with recurrence coefficients 
$(\beta_n^*)_n$ and $(\gamma_n^*)_n$  where 
\begin{equation}\label{beta*}
\beta_n^*=\beta_n+r_{1,n}-r_{1,n+1}, \quad n \ge 0,  
\end{equation}
\begin{equation}\label{gamma*}
\gamma_n^*=\gamma_n+r_{1,n} (\beta_{n-1}-\beta_n^*)+r_{2,n}-r_{2,n+1}, \quad n \ge 1,
\end{equation}
if and only if $\gamma_i^* \not= 0, \, i=1, \dots, N$ and the following formulas  hold:

\begin{equation}\label{condiciones P_{n-i}}
A_{i,n}=0, \quad n \ge i, \quad 2 \le i \le  N-1,\quad N \ge 3
\end{equation}
\begin{equation} \label{condicion P_{n-N}}
A_{N,n}=0, \quad n\ge N, \quad N \ge 2
\end{equation}
\begin{equation} \label{condicion P_{n-1-N}}
A_{N+1,n}=0 , \quad n\ge N+1, \quad N\ge 1
\end{equation} 
where
\begin{equation*} \label{A_n,i}
 A_{i,n}:=r_{i+1,n+1}-r_{i+1,n}+r_{i,n}(\beta_n^*-\beta_{n-i})+r_{i-1,n-1}\gamma_n^*- r_{i-1,n}\gamma_{n+1-i},
 \end{equation*}
\begin{equation*} \label{A_N,n}
A_{N,n}:=r_{N,n}(\beta_n^*-\beta_{n-N})+r_{N-1,n-1}\gamma_n^*- r_{N-1,n}\gamma_{n+1-N},
\end{equation*}
\begin{equation*}\label{A_N+1,n}
A_{N+1,n}:=r_{N,n-1}\gamma_n^*- r_{N,n}\gamma_{n-N}. 
\end{equation*}
\end{proposition}

\textbf{Proof.} We will characterize when the sequence $( R_n )_n$ is a MOPS, that is when it satisfies a three-term recurrence relation as
\begin{align} \label{recurrenciaRn}
R_{n+1}(x)&=(x-\beta_n^*)R_n(x)-\gamma_n^*R_{n-1}(x), \quad n\ge 0,
\end{align}
with  $\gamma_n^* \not=0, \quad n \ge 1$.

 Inserting formula (\ref{recurrenciaPn}) in
(\ref{relation Rn-Pn}) and applying (\ref{relation Rn-Pn}) to
$xP_n(x)$ and again (\ref{recurrenciaPn}) to $xP_{n-i}, i=1, \dots, N$,  and (\ref{relation Rn-Pn}) to $P_n(x)$ and next  to $P_{n-1}(x)$ we have for $n\ge 1$
\begin{align}\label{Rn+1-N3}
&R_{n+1}(x)=(x-\beta_n^*)R_n(x)-\gamma_n^*R_{n-1}(x)\\ \notag
&+ \sum_{i=2}^{N-1} [ r_{i+1,n+1}-r_{i+1,n}-r_{i,n}(\beta_{n-i}-\beta_n^*)-r_{i-1,n}\gamma_{n+1-i}+r_{i-1,n-1}\gamma_n^* ]P_{n-i}(x)\\ \notag
&+[r_{N,n}(\beta_n^*-\beta_{n-N})-r_{N-1,n}\gamma_{n+1-N}+r_{N-1,n-1}\gamma_n^*]P_{n-N}(x)\\ \notag
&+[-r_{N,n}\gamma_{n-N}+r_{N,n-1}\gamma_n^*]P_{n-1-N}(x),\ \notag
\end{align}
using (\ref{beta*}) and (\ref{gamma*}).

As a first consequence  for $n=1$ and any $N\ge 1$,  we realize that (\ref{recurrenciaRn}) is true with $\gamma_1^* \not=0$ if and only if $\gamma_1^* \not=0$.

Now, taking into account that the sequence  $( P_n )_n$ is a basis, we have that (\ref{recurrenciaRn})  with $\gamma_n^* \not=0$  for $ n \ge 2$ holds if and only if $\gamma_i^* \not=0, \, 2 \le i \le N$ and the  conditions 
(\ref{condiciones P_{n-i}}), (\ref{condicion P_{n-N}}) and (\ref{condicion P_{n-1-N}}) are satisfied. 
\hfill$\Box$

\medskip

\noindent {\bf Remark}. Observe that the condition (\ref{condicion P_{n-1-N}}) assures that $\gamma_n^* \not=0$ for $n \ge N+1$.

\medskip

In the following proposition, we change the hypothesis and assume now that  $(R_n)_n$  is a MOPS and characterize the orthogonality of the others polynomials which appear in the linear combination.

\begin{proposition}\label{caracterizacion ortogonalidad Qn} Let $(R_n)_n$ be a MOPS with  recurrence coefficients  $(\beta_n^*)_n$ and $(\gamma_n^*)_n$. Fixed $M \ge 1$, we define a sequence
$(Q_n)_n$ of monic polynomials recursively by 
\begin{equation} \label{relation Rn-Qn}
R_n(x):=Q_n(x)+\sum_{i=1}^M s_{i,n}Q_{n-i}(x)
\quad n\ge0,
\end{equation}
\noindent where  $s_{i,n}$ are complex numbers  such that $s_{i,n}=0$ for $i>n$ and $s_{M,n} \not=0$ for all $n \ge M$. 

Then $(Q_n)_n$ is a MOPS with recurrence coefficients $(\tilde{\beta}_n)_n$ and
$(\tilde{\gamma}_n)_n$ where 
\begin{equation}\label{betatilde}
\tilde{\beta}_n=\beta_n^*+s_{1,n+1} -s_{1,n}, \quad n \ge 0,
\end{equation}
\begin{equation}\label{gammatilde}
\tilde{\gamma}_n= \gamma_n^*+s_{1,n} (\beta_n^*-\tilde{\beta}_{n-1})+s_{2,n+1}-s_{2,n},  \quad n \ge 1,
\end{equation}
if and only if the following formulas hold: 
\begin{equation}\label{condiciones Q_{n-i}}
B_{i,n}=0, \quad n \ge i, \quad 2 \le i \le  M-1,\quad M \ge 3,
\end{equation}
\begin{equation} \label{condicion Q_{n-M}}
B_{M,n}=0, \quad n\ge M, \quad M \ge 2,
\end{equation}
\begin{equation}\label{condicion Q_{n-1-M}}
B_{M+1,n}=0,\quad n\ge M+1, \quad M \ge 1,
\end{equation}
\noindent {where}
\begin{equation} \label{B_n,i}
 B_{i,n}:=s_{i+1,n}-s_{i+1,n+1}+s_{i,n}(\tilde{\beta}_{n-i}-\beta_n^*)+s_{i-1,n}\tilde{\gamma}_{n+1-i}-s_{i-1,n-1}\gamma_n^*,
 \end{equation}
 \begin{equation} \label{B_M,n}
B_{M,n}:=s_{M,n}(\tilde{\beta}_{n-M}-\beta_n^*)+s_{M-1,n}\tilde{\gamma}_{n+1-M}-s_{M-1,n-1}\gamma_n^*,
\end{equation}
\begin{equation}\label{B_M+1,n}
B_{M+1,n}:=s_{M,n}\tilde{\gamma}_{n-M}-s_{M,n-1}\gamma_n^*.
\end{equation}
\end{proposition}

\textbf{Proof.}  Inserting the three-term recurrence relation satisfied by the polynomials $R_n$  and applying (\ref{relation Rn-Qn}) to
$xR_n(x)$,  $R_n(x)$ and $R_{n-1}(x)$, successively,  we have for $n\ge2$
\begin{align}\label{formula anterior empale}
Q_{n+1}(x)&=x\left[Q_n+\sum_{i=1}^M s_{i,n}Q_{n-i}(x)\right]-\beta_n^*\left[Q_n(x)+\sum_{i=1}^M s_{i,n}Q_{n-i}(x)\right]\\\notag
&-\gamma_n^*\left[Q_{n-1}(x)+\sum_{i=1}^M s_{i,n-1}Q_{n-1-i}(x)\right]-\sum_{i=1}^M s_{i,n+1}Q_{n+1-i}(x). \\\notag
\end{align}
Then, we can write 
\begin{align} \label{formula empalme}
&Q_{n+1}(x)=\notag\\
&xQ_n-(\beta_n^*+s_{1,n+1}-s_{1,n})Q_n(x)-(\gamma_n^*+s_{1,n}\beta_n^*+s_{2,n+1}-s_{2,n})Q_{n-1}(x)\notag\\ 
&+\sum_{i=1}^M s_{i,n} \left[xQ_{n-i}(x)-Q_{n+1-i}(x)  \right]+\sum_{i=3}^M s_{i,n}Q_{n+1-i}(x)\\
&-\beta_n^*\sum_{i=2}^M s_{i,n}Q_{n-i}(x)-\sum_{i=3}^M s_{i,n+1}Q_{n+1-i}(x)-\gamma_n^*\sum_{i=1}^M s_{i,n-1}Q_{n-1-i}(x).\notag
\end{align}

Now we suppose  that $(Q_n)_n$ is a MOPS with recurrence coefficients $(\tilde{\beta}_n)_n$ and
$(\tilde{\gamma}_n)_n$, that is the polynomials satisfy
\begin{equation} \label{recurrenciaQn}
Q_{n+1}(x)=(x-\tilde{\beta}_n)Q_n(x)-\tilde{\gamma}_nQ_{n-1}(x), \quad n\ge0,
\end{equation}
with  $\tilde{\gamma}_n \not=0$ for all $n\ge1$. Hence, applying this recurrence relation  to every  factor $[xQ_{n-i}(x)-Q_{n+1-i}(x)] $ in formula (\ref {formula empalme}) and using  (\ref{betatilde}) and (\ref{gammatilde}), it can be derived for $n\ge 2$
\begin{equation}\label{Qn+1-M3}
\sum_{i=2}^{M-1} B_{i,n}Q_{n-i}(x)+B_{M,n}Q_{n-M}(x)+B_{M+1,n}Q_{n-1-M}(x)=0.
\end{equation}

Then, the conditions (\ref{condiciones Q_{n-i}}),  (\ref{condicion Q_{n-M}})  and  (\ref{condicion Q_{n-1-M}}) hold, for $n\ge2$.

\medskip

In order to proof the reverse, we observe that formula (\ref{formula empalme}) and  conditions (\ref{condiciones Q_{n-i}}),  (\ref{condicion Q_{n-M}})  and  (\ref{condicion Q_{n-1-M}}) yield

\begin{align} \label{FRQn+1--FRQn+1-i}
&Q_{n+1}(x)-(x-\tilde{\beta}_n)Q_n(x)+\tilde{\gamma}_nQ_{n-1}(x)\\\notag
&=-\sum_{i=1}^M s_{i,n}\left[Q_{n+1-i}(x)-(x-\tilde{\beta}_{n-i})Q_{n-i}(x)+\tilde{\gamma}_{n-i}Q_{n-1-i}(x) \right],\\\notag
\end{align}
for $n\ge2$.

Since formula (\ref{recurrenciaQn}) is true for $n=1$ and obviously for $n=0$, we have that it is  also true for $n=2$. Hence by using the induction method  we can deduce that  the recurrence relation (\ref{recurrenciaQn}) holds for  all $n \ge 0$. 

To conclude it suffices to observe that   $\tilde{\gamma}_n\not=0$ for $n \ge 1$,  from  (\ref{condicion Q_{n-1-M}}). \hfill$\Box$

\medskip

In the following theorem, we link the results studied in Proposition \ref{caracterizacion ortogonalidad Rn} and Proposition \ref{caracterizacion ortogonalidad Qn} and thus we get to solve the inverse problem in a general case without many calculations. 

\begin{theorem}\label{caracterizacion ortogonalidad caso general} 
Let  $( P_n )_n$ be a MOPS with respect to a regular functional $\bf u$ with $(\beta_n)_n$ and $(\gamma_n)_n$ the corresponding sequences of recurrence coefficients. We define recursively a sequence  $( Q_n )_n$ of monic polynomials by formula (\ref{relation N-M}), i.e.
$$P_n(x)+\sum_{i=1}^N r_{i,n}P_{n-i}(x)=Q_n(x)+\sum_{i=1}^M s_{i,n}Q_{n-i}(x), \quad n \ge 0, $$
where $ (r_{i,n})_{i=1}^N$ and $(s_{i,n})_{i=1}^M$ are  complex numbers and such that
 the conditions  $\text{det}\, {\bf A} \not=0$ and  $r_{N,n}\,s_{M,n}\not=0$  for all $ n \ge N+M$, are satisfied.
 
 Then $( Q_n )_n$ is a MOPS with recurrence coefficients $(\tilde{\beta}_n)_n$ and $(\tilde{\gamma}_n)_n, $  where 
 \begin{equation}\label{betatilde-beta}
 \tilde{\beta}_n+s_{1,n}-s_{1,n+1}=\beta_n+r_{1,n}-r_{1,n+1}, \quad n \ge0,
\end{equation}
\begin{align}\label{gammatilde-gamma}
& \tilde{\gamma}_n+ s_{1,n} \left(\tilde{\beta}_{n-1}-\tilde{\beta}_n+s_{1,n+1}-s_{1,n} \right)+s_{2,n}-s_{2,n+1}\\ \notag
&= {\gamma}_n+ r_{1,n} \left({\beta}_{n-1}-{\beta}_n+r_{1,n+1}-r_{1,n} \right)+r_{2,n}-r_{2,n+1}, \quad n\ge1.
\end{align}
 if and only if the polynomials $Q_n$ satisfy the three-term recurrence relation  with $\tilde{\gamma}_n\not=0$ for $n=1,2, \dots, N$ and the equations  (\ref{condiciones P_{n-i}})--(\ref{condicion P_{n-1-N}}) and (\ref{condiciones Q_{n-i}})--(\ref{condicion Q_{n-1-M}})  for $n \ge N+M+1$ hold.
\end{theorem}

\textbf{Proof.} Inserting formula (\ref{recurrenciaPn}) in (\ref{relation N-M}) and applying  (\ref{relation N-M}) to $xP_n(x)$ and again (\ref{recurrenciaPn}) to $xP_{n-i}(x)$ for $i=1,2, \dots, N$ we get, for $n\ge1$,
\begin{align*}
Q_{n+1}(x)&=-\sum_{i=1}^Ms_{i,n+1}Q_{n+1-i}(x)+x\left[Q_n(x)+\sum_{i=1}^Ms_{i,n}Q_{n-i}(x) \right]\\
&-\beta_nP_n(x)-\gamma_nP_{n-1}(x)+\sum_{i=1}^Nr_{i,n+1}P_{n+1-i}(x)\\
&-\sum_{i=1}^Nr_{i,n} \left[P_{n+1-i}(x)+\beta_{n-i}P_{n-i}(x)+\gamma_{n-i}P_{n-1-i}(x) \right].
\end{align*}

Now, in the above expression, we apply (\ref{relation N-M}) to $P_n(x)$ and we rearrange the formula.  Next, we do the same for $P_{n-1}(x)$. Hence,  using  the auxiliary coefficients (\ref{beta*}) and (\ref{gamma*}), we can write the above formula in the following way

\begin{align}\label{empalme Pn-Qn}
Q_{n+1}(x)&=-\sum_{i=1}^Ms_{i,n+1}Q_{n+1-i}(x)+x \left[Q_n(x)+\sum_{i=1}^Ms_{i,n}Q_{n-i}(x) \right]\\ \notag
&-\beta_n^* \left[Q_n(x)+\sum_{i=1}^Ms_{i,n}Q_{n-i}(x) \right]\\ \notag
&-\gamma_n^* \left[Q_{n-1}(x)+\sum_{i=1}^Ms_{i,n-1}Q_{n-1-i}(x)) \right]\\
&+\sum_{i=2}^{N-1}A_{i,n} P_{n-i}(x)+A_{N,n}P_{n-N}(x)+A_{N+1,n}P_{n-(N+1)}(x).\notag
\end{align}

Let  $(Q_n)_n$ be a MOPS with recurrence coefficients $(\tilde{\beta}_n)_n$ and
$(\tilde{\gamma}_n)_n$. Then, applying (\ref{recurrenciaQn}) to every  factor $xQ_{n-i}(x)$ for $i=1, \dots, M$ in formula (\ref{empalme Pn-Qn})  and using  (\ref{betatilde}), (\ref{gammatilde}) and  (\ref{B_n,i})--(\ref{B_M+1,n}
), we get
\begin{align}\label{relaci—n Ai-Bi}
0&=\sum_{i=2}^{M-1}B_{i,n}Q_{n-i}(x)+B_{M,n}Q_{n-M}(x)+B_{M+1,n}Q_{n-(M+1)}(x)\\\notag
+&\sum_{i=2}^{N-1}A_{i,n}P_{n-i}(x)+A_{N,n}P_{n-N}(x)+A_{N+1,n}P_{n-(N+1)}(x).
\end{align}

Now applying the functional $\overline{P}_{M-1}{\bf u}$ to the above equation and taking into account the orthogonality of the polynomials $P_n$ with respect to $\bf u$,  we obtain 
$$\langle \overline{P}_{M-1}{\bf u}, \sum_{i=2}^{M-1}B_{i,n} Q_{n-i}+B_{M,n}Q_{n-M}+B_{M+1,n}Q_{n-(M+1)} \rangle=0, \quad n \ge N+M+1,$$
$$\langle \overline{P}_{M-2}{\bf u}, \sum_{i=2}^{M-1}B_{i,n} Q_{n-i}+B_{M,n}Q_{n-M}+B_{M+1,n}Q_{n-(M+1)}  \rangle=0, \quad n \ge N+M,$$
and successively until
$$\langle \overline{P}_0{\bf u}, \sum_{i=2}^{M-1}B_{i,n} Q_{n-i}+B_{M,n}Q_{n-M}+B_{M+1,n}Q_{n-(M+1)}  \rangle=0, \quad n \ge N+2.$$
Thus, for $n \ge N+M+1$ we have
 \begin{equation}\label{sumas Bi}
 \sum_{i=2}^{M-1}B_{i,n} \langle  \overline{P}_{j}\, {\bf u}, Q_{n-i} \rangle+B_{M,n} \langle  \overline{P}_{j}\, {\bf u},Q_{n-M} \rangle+
B_{M+1,n} \langle  \overline{P}_{j}\, {\bf u},Q_{n-(M+1)} \rangle=0
\end{equation}  
for $ j=0,1, \dots, M-1.$

Fixed an positive integer $n$ ( $n \ge N+M+1$), we consider the homogeneous system (\ref{sumas Bi}) of $M$ equations and $M$ unknowns $B_{i,n}, \, i=2, \dots, M+1,$ whose associated matrix is ${\widetilde{\bf B}}^T_{n-2}$ that is the transpose of the matrix $\widetilde{\bf B}_{n-2}$.  Applying  the hypotheses 
$\text{det} {\bf A} \not=0$, $\,s_{M,n}\not=0$  for all $ n \ge N+M$ and $(b)$ of the Lemma \ref {relation entre Bn}, we get $\text {det}\,{\widetilde{\bf B}}^T_{n-2}\not=0$, for all $n \ge N+M+1$. Therefore the system has a unique solution that is the trivial solution 
$$B_{i,n}=0,\quad  i=2, \dots, M+1,$$ so the equations (\ref{condiciones Q_{n-i}})--(\ref{condicion Q_{n-1-M}}) hold for $n \ge N+M+1.$   The other equations (\ref{condiciones P_{n-i}})--(\ref{condicion P_{n-1-N}}) for $n \ge N+M+1$ are also fulfilled simply noticing that from the equation (\ref{relaci—n Ai-Bi})  we obtain
$$\sum_{i=2}^{N-1}A_{i,n} P_{n-i}(x)+A_{N,n}P_{n-N}(x)+A_{N+1,n}P_{n-(N+1)}(x)=0$$ and the sequence of polynomials $(P_n)_n$ is a basis.

Conversely, notice that the equations (\ref{condiciones P_{n-i}})--(\ref{condicion P_{n-1-N}})  in formula (\ref{empalme Pn-Qn})  yield
\begin{align*}
Q_{n+1}(x)&=-\sum_{i=1}^Ms_{i,n+1}Q_{n+1-i}(x)+x \left[Q_n(x)+\sum_{i=1}^Ms_{i,n}Q_{n-i}(x)  \right]\\
&-\beta_n^* \left[Q_n(x)+\sum_{i=1}^Ms_{i,n}Q_{n-i}(x)) \right] 
-\gamma_n^* \left[Q_{n-1}(x)+\sum_{i=1}^Ms_{i,n-1}Q_{n-1-i}(x) \right]
\end{align*}
 for $n \ge N+M+1.$  We observe that this is  formula (\ref {formula anterior empale}) in the Proposition \ref{caracterizacion ortogonalidad Qn} which can be rewritten as (\ref{formula empalme}).  Then taking into account the equations (\ref{condiciones Q_{n-i}})--(\ref{condicion Q_{n-1-M}})  we get 
(\ref{FRQn+1--FRQn+1-i}) for $n \ge N+M+1.$ Thus, since by hypotheses the polynomials $Q_n$ satisfy the three-term recurrence relation  with  $\tilde{\gamma}_n\not=0$ for $n=1,2, \dots, N$ we obtain that these polynomials satisfy the same  three-term recurrence relation for all $n$.  Besides  from (\ref{condicion P_{n-1-N}}) and (\ref{condicion Q_{n-1-M}}) and  we get that $\tilde{\gamma}_n\not=0$ for all $n \ge N+1 $  and therefore the sequence $( Q_n )_n$ is a MOPS.  \hfill$\Box$

\medskip

\noindent {\bf Remark}.  Observe that to prove the first part of the Theorem, i.e. with the assumption that the sequence  $(Q_n)_n$ is orthogonal,  it is enough to require only the initial conditions (\ref {initial conditions}), taking into account (a) of Theorem \ref{unicidad Psi-Phi}. However for the converse we need that the more extensive condition   $r_{N,n}\,s_{M,n}\not=0$ for all $ n \ge N+M$  hold.

\end{document}